\documentclass[12pt]{article}
\newtheorem{theorem}{Theorem} 
\newtheorem{lemma}[theorem]{Lemma} 
\newtheorem{proposition}[theorem]{Proposition}

\def\bthm#1.#2 #3\ethm{ 
\begin{\ifx#1ttheorem\fi%
\ifx#1llemma\fi%
\ifx#1ccorollary\fi%
\ifx#1pproposition\fi%
\ifx#1ddefinition\fi} \label{#1.#2} 
#3 \end{\ifx#1ttheorem\fi%
\ifx#1llemma\fi%
\ifx#1ccorollary\fi%
\ifx#1pproposition\fi%
\ifx#1ddefinition\fi}} 
 
\def\t#1/{theorem~\ref{t#1}}   \def\T#1/{Theorem~\ref{t#1}} 
\def\c#1/{corollary~\ref{c#1}}   \def\C#1/{Corollary~\ref{c#1}} 
\def\l#1/{lemma~\ref{l#1}}        
\def\s#1/{section~\ref{s#1}} 
\def\e#1/{(\ref{e#1})} 
\def\d#1/{definition~\ref{d#1}} 
\def\f#1/{figure~\ref{f#1}} 
 
\def\Label #1 {\label{#1}}

\def\norm#1.#2.{\lVert#1\rVert_{#2}} 
\def\Norm#1.#2.{\bigl\lVert#1\bigr\rVert_{#2}} 
\def\NOrm#1.#2.{\Bigl\lVert#1\Bigr\rVert_{#2}} 
\def\NORm#1.#2.{\biggl\lVert#1\biggr\rVert_{#2}} 
\def\NORM#1.#2.{\Biggl\lVert#1\Biggr\rVert_{#2}} 
 
\def\ip#1,#2.{\langle #1,#2\rangle} 
\def\Ip#1,#2.{\bigl\langle#1,#2\bigr\rangle} 
\def\IP#1,#2.{\Bigl\langle#1,#2\Bigr\rangle}

\def\ZZ{{\bf Z}}

\def\mid{\,:\,} 
 \newcommand{\qed}{{\vrule height 1ex width 1ex depth -.1ex}}
 
\def\md#1#2\emd{\ifx0#1 
\begin{equation*} #2 \end{equation*}\fi  
\ifx1#1\begin{equation}#2\end{equation}\fi   
\ifx2#1\begin{align*}#2\end{align*}\fi   
\ifx3#1\begin{align}#2\end{align}\fi    
\ifx4#1\begin{gather*}#2\end{gather*}\fi  
\ifx5#1\begin{gather}#2\end{gather}\fi   
\ifx6#1\begin{multline*}#2\end{multline*}\fi  
\ifx7#1\begin{multline}#2\end{mutline}\fi  
}

\begin{document}

\title{  On  Sets of Integers 
\\ Not Containing Long Arithmetic Progressions }

 \author{
 Izabella {\L}aba\thanks{Supported in part by NSERC}
 \\ University of British Columbia
   \and
 Michael T. Lacey\thanks{Supported by an NSF grant, DMS--9706884.}\\
 Georgia Institute of Technology 
 }

\maketitle


\section{The Main Result}
\label{sec.1}

Let $r(k,N)$ be the maximal cardinality of a subset $A$ of $\{1,2,\ldots,N\}$ 
which does not contain an arithmetic progression of length $k$.  That is, $A$ 
does not contain a subset of the form $\{x+jy\mid 0\le{}j<k\}$, where $x,y      $ 
are integers with $y\not=0$. 

Erd\"os and Turan \cite{et}
initiated the study of these   quantities  in 1936.  In particular they 
conjectured that $r(k,N)=o(N)$ for all $k$, that is every set of integers of 
positive asymptotic density contains arbitrarily long arithmetic progressions.  
In 1953, Roth \cite{R1} showed that $r(3,N)=o(N)$. 
The Erd\"os--Turan conjecture was verified by Szemer\'edi \cite{Sz1,Sz2}, 
a result with a very broad influence.  Subsequently, rather different proofs
of Szemer\'edi's theorem were given by Furstenberg \cite{Fu} and Gowers \cite{g1,g2}.
Gowers's proof provides, for the first time, upper bounds on $r(k,N)$ given by 
a bounded tower of exponentials.
 An intriguing question of Erd\"os asks if $r(3,N)\le{}CN/(\log N)^{1+\delta}$ 
for some positive $\delta$.
 Bourgain's article \cite{B} contains the best 
current upper bound of $CN\sqrt{\frac{\log\log N}{\log N}}$ on $r(3,N)$ .  

In this article we are interested in the converse question of finding large
subsets of $\{1,\dots,N\}$ which do not contain arithmetic progressions. 
 Behrend, in 1946, 
  \cite{b} (building on earlier work of Salem and Spencer \cite{ss}) 
  considered three term arithmetic progressions, and showed 
 that $r(3,N)\ge{}N \,\exp(-C\sqrt{\log N})$.  
The purpose of this paper is to show that if one considers longer
 arithmetic progressions then Behrend's estimate can be further improved as follows.

\begin{theorem}\label{t.antiroth}
There is a constant $C>0$ so that for all $n>k\geq 1$, 
\begin{equation}\label{e.main}
 r(1+2^k,N)\ge{}N\exp(-C (\log N )^{1/(k+1)}).
\end{equation}
\end{theorem}

\section{The Proof}
\label{sec.2} 

Our argument builds upon the methods of Salem and Spencer \cite{ss} and of 
Behrend \cite{b}.  It will be convenient to consider the set $I=\ZZ\cap(-
\frac{N-1}{2},\frac{N-1}{2}]$ instead of $\{1,2,\dots,N\}$.  
First, we may assume that $N=n^d$ for suitably chosen integers $n$ and $d$, 
with $n$ much smaller than $N$ and divisible by a constant $c_0$ (independent
of $N,n$) to be chosen later.  Indeed, 
at the cost of a slightly larger constant in our theorems we can always increase 
$n$ to one of these values.  Similarly, we shall take 
fractional powers and logarithms of large integers and tacitly assume that the output is 
also an integer.  In fact the argument requires  the integer parts of these 
quantities, but to minimize notation we do not explicitly invoke the integer 
part function.  

Second, with $N=n^d$, consider the expansion of each $x\in I$ in base 
$n$, defined as follows.  For any $x\in I$ we define its coordinate vector 
$v_x=(x_0,\dots,x_{d-1})\in\ZZ^d$, where $x_i$ are uniquely determined by the
conditions
\begin{equation}\label{e.digits}
x=\sum_{i=0}^{d-1} x_i n^i,\ -\frac{n-1}{2}< x_i\leq \frac{n-1}{2}.
\end{equation}
Note that, unlike in Behrend's argument, the ``digits" $x_i$ are not
required to be non-negative.  
Denote also the ``norm" of $x\in I$ as 
\[
\|x\|^2=\|v_x\|^2=\sum_{i=0}^{d-1}x_i^2,
\]
with $x_i$ defined by (\ref{e.digits}).

An important observation of Salem and Spencer \cite{ss} was that if we only 
consider the set $Q_0$ of numbers $x\in\{0,1,\dots,N-1\}$ with digits 
$0\le{}x_i<cn$, where $c$ is sufficiently small\footnote{Salem and Spencer 
considered expansions with non-negative digits $0\leq x_i\leq n-1$, in which
case it suffices to take $c=1/2$}, then 
addition of numbers is equivalent to vector addition in the corresponding subset
of $\ZZ^d$, i.e. for any $x,y,z\in Q_0$ we have $x+y=z$ if and only if $v_x+v_y
=v_z$.  Thus an arithmetic progression $x,x+y,x+2y,\dots$ in $Q_0$ corresponds to 
vectors $v_x, v_{x+y}, v_{x+2y},\dots$ on a straight line in $\ZZ^d$.  

We shall rely on variants of this observation.  More precisely, we define
\begin{equation}\label{e.Q}
Q=\{x\in\ZZ:\ x=\sum_{i=0}^{d-1} x_i n^i,\ 
-q\leq x_i\leq q\},
\end{equation}
where $q=n/c_0$ and $c_0$ is a large constant independent of $N,n$ to be
chosen later. We will also denote for $r\in\ZZ$:
\[
rQ=\{x\in\ZZ:\ x=\sum_{i=0}^{d-1} x_i n^i,\ 
-rq\leq x_i\leq rq\}.
\]
Then linear combinations of numbers in $rQ$ with small enough integer 
coefficients correspond to linear combinations of their coordinate vectors:
\begin{equation}\label{e.lincomb}
v_{\sum a_kx^{(k)}}=\sum a_kv_{x^{(k)}}
\hbox{ if }a_k,r_k\in\ZZ,\ x^{(k)}\in r_kQ,\ \sum r_k|a_k|<c_0/3.
\end{equation}

Our proof consists of two distinct parts, both similar in spirit to Behrend's
argument \cite{b}.  The latter  relies on the geometrical fact 
that a straight line can intersect a sphere $\|v_x\|^2=r$ in $\ZZ^d$ in at most two points, 
so that the set $\{x\in Q:\ \|x\|^2=r\}$ cannot contain a three-term arithmetic 
progression. One then uses pigeonholing to choose a sphere containing a large 
number of points in $Q$.  

Our intermediate results can be stated in terms of quantities closely related
to those of Erd\"os and Turan.  Namely, define
$r_m(k,N)$ to be the maximal cardinality of a subset $A\subset \{0,1,\ldots,N-1\}$ 
which does not contain a further subset of the form 
\begin{equation}\label{e.excludedsets}  
 \Bigl\{ x+\sum_{i=1}^{m} a_i j^i\mid 0\le{}j<k-1 \Bigr\},
\end{equation}
for any integers $x$ and $a_i$ such that at least one of the $a_i$ is non-zero.
(In particular, $r_1(k,N)=r(k,N)$ and $r_m(k,n)$ decreases with $m$.)
 Observe that a set of the form 
(\ref{e.excludedsets}) with $u\geq 2$ may contain less than $k$ 
distinct integers, as the same summand may arise from  more than one
value of $j$.  Note further that the $a_i$ need not belong to $A$. 
Finally, while this is defined as a property of the initial interval of 
integers $\{0,\ldots,N-1\}$,
it depends only on the length of the interval of integers in question.
 
The estimates we will need are the following.

\begin{proposition}\label{p.parabola}
We have
\begin{equation}\label{e.r2p}  
r_m(2m+1,N)\ge  N\exp(-C\sqrt{\log N}),
\end{equation}
where $C$ is an absolute constant depending only on $m$.
\end{proposition}

\begin{proposition}\label{prop.induct}
Assume that $N=n^d$, and let $k\geq m+1$. Then 
\begin{equation}\label{e.induct}
  r_m(k,N)\geq N\frac{r_{2m}(k,n^2d)}{c^d\,n^2d},
\end{equation}
where the constant $c>0$ depends only on $m$ and $k$.
\end{proposition}

Proposition \ref{p.parabola} is proved by essentially repeating Behrend's 
argument with straight lines
replaced by curves of higher order; the main point is that a non-constant
polynomial of degree $2m$ can have at most $2m$ roots.
Proposition \ref{prop.induct} will allow us to carry out the inductive argument.
Instead of just one sphere as in Behrend's argument, the set $A$ which provides
the lower bound in (\ref{e.induct}) will be a union of concentric spheres of radii
$\sqrt{r}$, $r\in R$.  We will argue that if $A$ contains a subset $\{x^{(j)}\}$ as in
(\ref{e.excludedsets}), then the squared norms $\|x^{(j)}\|^2$ are
as in (\ref{e.excludedsets}) with $m$ replaced by $2m$.  Proposition
\ref{prop.induct} will follow upon choosing a set $R$ of cardinality 
$r_{2m}(k,dn^2)$ which cannot contain such a subset, and optimizing over $n$ and $d$.

We will use $C, c, c_i,$ etc. to denote absolute constants which may depend on $m$ and
may change from line to line but are always independent of $N$, $n$, $d$.

\section{Proof of Proposition \ref{p.parabola}}
\label{sec.3}

Our goal in this section is to find a set $R\subset \{0,1,\dots,N-1\}$ of large 
cardinality such that $R$ does not contain all of the integers 
\begin{equation}\label{e.5par}
\sum_{i=0}^m a_ij^i:\ j=0,1,\dots,2m
\end{equation}
for any $a_0,\dots,a_m\in\ZZ$ with $a_i\neq 0$ for at least one $i>0$.
We will use the notation of Section \ref{sec.2}.  In particular, we will
replace the set $\{0,1,\dots,N-1\}$ by $I$, and
assume that $N=n^d$ for some $1\ll d\ll N$ and $1\ll n
\ll N$ (eventually we will let $d\sim\sqrt{\log N}$).  The set $R$  
will be a subset of the set $Q$ defined in (\ref{e.Q}).  

\begin{lemma}\label{l.smallabc}
Suppose that $2m+1$ numbers $x^{(j)}$ in $Q$ satisfy
\begin{equation}\label{e.five1}
x^{(j)}=\sum_{i=0}^m a_ij^i,\ j=0,1,\dots,2m,
\end{equation}
for some integers $a_0,\dots,a_m$.  Denote by $D$ the Vandermonde 
determinant $D=D_m=|J_m|$, where $J_m=(j^i)_{i,j=1}^m$. 
Then there is a constant $c$, depending only on $m$, such that 
\begin{equation}\label{e.five3}
D a_i\in cQ,\ i=0,\dots,m.
\end{equation}
Furthermore, if the constant $c_0$ in the definition of $Q$ was chosen 
large enough, then we have for any such numbers
\begin{equation}\label{e.five2}
D v_{x^{(j)}}=\sum_{i=0}^m j^i\,v_{D a_i},\ j=0,1,\dots,2m.
\end{equation}
\end{lemma}

{\it Proof.} 
We consider the first $m+1$ equations in (\ref{e.five1}) as a system of 
linear equations with unknowns $a_0,\dots,a_m$. By Cramer's formula,
$Da_i$ are linear combinations of $x^{(j)}$ with integer coefficients
bounded by a constant depending only on $m$.  This implies 
(\ref{e.five3}). Now (\ref{e.five2}) follows from (\ref{e.five3}), 
(\ref{e.five1}) and (\ref{e.lincomb}). 
\qed

\bigskip
We are now in a position to run Behrend's argument.  Let
\[
S_r=\{x\in Q:\ \|x\|^2=r\},
\]
where $\|x\|^2=\|v_x\|^2=\sum_{i=0}^d |x_i|^2$.  We will prove that
no $S_r$ may contain $2m+1$ points as in (\ref{e.five1}). Indeed, suppose to
the contrary that $x^{(j)}$, $j=0,1,\dots,2m$, satisfy (\ref{e.five1}) and
$\|x^{(j)}\|^2=r$.  By Lemma \ref{l.smallabc}, we have
\[  P(j):={}
\|x^{(j)}\|^2=\sum_{k=0}^{d-1}\Big(\sum_{i=0}^m \frac{(Da_i)_k}{D}j^i \Big)^2.
\]
But then $P(j)$ is a polynomial of degree $2m$ in $j$, equal to $r$ for 
$j=0,1,\dots,2m$.  This is not possible unless $P(j)$ is constant, in
which case we must have $(Da_i)_k=0$ for all $0\leq k\leq d-1$ and all
$1\leq i\leq m$.  By Lemma \ref{l.smallabc} again, it follows that $a_i=0$
for all $1\leq i\leq m$.

Finally, we use a pigeonholing argument to find a set $S_r$ of large cardinality.
Following Behrend \cite{b}, we set $d=\sqrt{\log N}$ and $n=N^{1/d}$, so that
$q=N^{1/d}/1000$. Since $Q$ has cardinality $(2q)^d$ and $Q=\bigcup_{r=0}^{dq^2}
S_r$, there is at least one $r$ for which
\[                                 
N^{-1}\# S_r\geq{} (d 500^dq^2)^{-1}\geq{}C_1\exp    (-C_1d)N^{-2/d}
{}\geq{}C_1\exp\Big(-C_1\Big(d+\frac{\log N}d\Big)\Big).
\]
Taking $d=\sqrt{\log N}$ proves the proposition.

\section{Proof of Proposition \ref{prop.induct}} 
\label{sec.4}

We continue to use the notation of Section \ref{sec.2}: we assume that
$N=n^d$ with $n,d\ll N$,
and define $q$, $Q$, $v_x$, $\|x\|$, $D$, etc. as before.  We also define
\[
(x,y)=\sum_{i=0}^{d-1} x_iy_i
\]
for $x,y\in \frac{c_0}{3} Q$.  

Let $R\subset \{0,1,\dots,D^2dq^2-1\}$ be a set of cardinality $r_{2m}(k,D^2dq^2)$ which
does not contain all of the integers
\begin{equation}\label{e.radii}
  y^{(j)}=\sum_{i=0}^{2m}a_ij^i,\ j=0,1,\dots,k-1,
\end{equation}
for any $a_0,\dots,a_{2m}\in\ZZ$.  
Observe that any translate $R+s:=\{r+s:\ r\in R\}$, $s\in\ZZ$, of $R$ has
the same cardinality as $R$ and cannot contain $k$ integers as in (\ref{e.radii}).
Let $X=2q\sum_{i=0}^{d-1}n^i\in 2Q$ and
$S:=\{0,1,\dots,9 D^2dq^2\}$.  For $s\in S$,  define
\[
A_s=\{x\in Q:\ D^2\|x-X\|^2\in R+s\}.
\]
We claim that no $A_s$ can contain $k$ integers
\begin{equation}\label{e.lower}
  x^{(j)}=\sum_{i=0}^{m} b_ij^i,\ j=0,1,\dots,k-1.
\end{equation}
Indeed, suppose to the contrary that $A_s$ does contains such $k$ integers.  
As in Lemma \ref{l.smallabc}, we prove that 
\[
D v_{x^{(j)}-X}=\sum_{i=0}^{m}j^i\,v_{Db_i}-Dv_{X},\ j=0,\dots,k-1,
\]
provided that $c_0$ was chosen large enough.  Hence
\[
D^2\|x^{(j)}-X\|^2=\sum_{k=0}^{d-1}\Big( \sum_{i=0}^{m} j^i(Db_i)_k-2Dq\Big)^2
\]
are as in (\ref{e.radii}).  But this is impossible by the choice of $R$. 

\medskip 

A pigeonholing argument shows that there is an $A_s$ with large cardinality. 
For any $x\in Q$ we have $q\leq (X-x)_i\leq 3q$ for each $i$, 
 hence $D^2dq^2\leq D^2\|x-X\|^2\leq 
9D^2dq^2$. Hence for any $x\in Q$ and $r\in R$ we have
\[
1\leq \|x-X\|^2-r\leq 9D^2dq^2,
\]
and in particular there is a $s\in S$ such that $D^2\|x-X\|^2=r+s$.  It follows that
for each $x\in Q$ there are at least $\# R$ values of $s$ such that $x\in A_s$.
Hence
\[
\sum_{s\in S}\# A_s\geq \# R\cdot \# Q.
\]
In particular, there is an $s\in S$ such that
\[
\# A_s\geq \frac{\# R\cdot \# Q}{\# S}\geq C\frac{n^d }{ 1000^{ d}}\cdot
\frac{r_m(k,D^2dn^2)}{D^2n^2d},
\]
which yields (\ref{e.induct}).

\section{Proof of Theorem \ref{t.antiroth}} 
\label{sec.5}

We will prove that for all $1\leq k\ll \log N$ and all
$1\leq l\leq k$,
\begin{equation}\label{e.kl}
  r_{2^{k-l}}(1+2^k,N)\geq N\exp\big(-c(\log N)^{\frac{1}{l+1}}\big).
\end{equation}
In particular, taking $l=k$ we obtain (\ref{e.main}).  Here and below,
the constants $c,c',c''$ may depend on $m,k,l$, but not on $N$.

The proof of (\ref{e.kl}) is by induction in $l$.  The case $l=1$ is
(\ref{e.r2p}).  Suppose now that (\ref{e.kl}) holds for $l$, and set
$N=n^d$, $d\sim (\log N)^{1/(l+2)}$. Then by (\ref{e.induct}) we have
\[
r_{2^{k-l-1}}(1+2^k,N)\geq N \frac{r_{2^{k-l}}(1+2^k,n^2d)}{c^d\,n^2d}
\geq Nc^{-d}\exp\big(-c'(\log (n^2d))^{\frac{1}{l+1}}\big)
\]
\[
\geq N\exp\big(-c''(\log N)^{\frac{1}{l+2}}\big),
\]
which is (\ref{e.kl}) for $l+1$.

\bigskip 
{\parindent=0pt\baselineskip=12pt\obeylines 
Izabella {\L}aba
Department of Mathematics
University of British Columbia 
Vancouver, B.C. V6T 1Z2, Canada 
\smallskip
\tt ilaba@math.ubc.ca
\tt http://www.math.ubc.ca/\~{}ilaba}

\bigskip
{\parindent=0pt\baselineskip=12pt\obeylines 
Michael T. Lacey 
School of Mathematics 
Georgia Institute of Technology 
Atlanta, GA 30332, U.S.A.
\smallskip
\tt lacey@math.gatech.edu
\tt http://www.math.gatech.edu/\~{}lacey }


\begin{thebibliography}{9} 
 
\bibitem{b} F.A. Behrend,  On sets of integers which contain no three terms in 
arithmetic progression,  Proc. Nat. Acad. Sci.  32 (1946), 331-332.

\bibitem{B} J.  Bourgain, On triples in arithmetic progression, Geom. Func. Anal. 
9 (1999) 968---984.

\bibitem{et}  P. Erd\"os and P. Turan,  On some sequences of integers, 
{     J. London Math. Soc. } {    11} (1936), 261---264. 

\bibitem{Fu} H. Furstenberg, Ergodic behaviour of diagonal measures and
  a theorem of Szemer\'edi on arithmetic progressions, J. Analyse Math. 31 (1977),
   204---256. 

\bibitem{g1}  W.T. Gowers, A new proof of Szemer\'edi's theorem for arithmetic 
progressions of length four, Geom. Func. Anal. 8 (1998), 529---551.


\bibitem{g2} W.T. Gowers,  A new proof of Szemer\'edi's theorem, 
            Geom. Func. Anal. 11 (2001), 465--588.
     
\bibitem{H} D.R. Heath-Brown, Integer sets containing no arithmetic 
progressions,     J. London Math. Soc. (2) 35 (1987), 385---394. 
    
\bibitem{R1} K.F. Roth, On certain sets of integers, J. London Math. Soc. 28 
(1953), 245---252. 
    
\bibitem{R2} K.F. Roth, Irregularities of sequences relative to 
arithmetic progressions, IV, Period. Math. Hungar. 2 (1972), 301---326. 
    
\bibitem{ss}  R. Salem and  D.C. Spencer, On sets of integers which contain no three terms in 
arithmetic progression,  Proc. Nat. Acad. Sci.  32 (1942),  561---563.


 
\bibitem{Sz1} E. Szemer\'edi, On sets of integers containing no four
elements in arithmetic progression, Acta Math. Acad. Sci. Hungar. 20 
(1969), 89---104. 
      
\bibitem{Sz2} E. Szemer\'edi, On sets of integers containing no k elements 
in arithmetic progression, Acta Arith. 27 (1975), 299---345. 
           
\end{thebibliography}
\end{document}